\documentstyle[theorem,amstex,amssymb]{article}

\newtheorem{proposition}{Proposition}[section]
\newtheorem{lemma}[proposition]{Lemma}
\newtheorem{theorem}[proposition]{Theorem}

{\theorembodyfont{\rm}
\newtheorem{definition}[proposition]{Definition}

\def\pmc#1{\setbox0=\hbox{#1}
    \kern-.1em\copy0\kern-\wd0
    \kern.1em\copy0\kern-\wd0}

\def\Wb{{\bar{W}}}
\def\Ub{{\bar{U}}}
\def\Vb{{\bar{V}}}
\def\WXb{{\bar{W}_X}}

\def\itheorem#1#2{\newtheorem{#1}[proposition]{#2}}

\itheorem{stacked}{Stacked Bases Theorem}

\itheorem{remark}{Remark}

\def\proof{\paragraph{Proof. }}
\def\qed{\hglue 0pt plus1fill $\Box$\\ \vspace{2ex}}

\def\Z{{\mathbb{Z}}}

\def\N{{\mathbb{N}}}

\begin{document}
\title{The failure of the uncountable non-commutative Specker Phenomenon
\author{By\\[0.2cm]Saharon Shelah \thanks{Number of publication 729.
Supported by project No. G-0545-173,06/97 of the {\em
German-Israeli Foundation for Scientific Research \& Development}}\quad
and Lutz Str\"ungmann\thanks{Supported by the Graduiertenkolleg {\em
Theoretische und Experimentelle Methoden der Reinen Mathematik} of Essen
University. }}
\makeatother
}
\date{}
\maketitle

\begin{abstract}
Higman proved in 1952 that every free group is non-commutatively
slender, this is to say that if $G$ is a free group and $h$ is a
homomorphism from the countable complete free product
$\pmc{$\times$}\, \, \;_{\omega} \Z$ to $G$, then there exists a
finite subset $F \subseteq \omega$ and a homomorphism $\bar{h} :
*_{i \in F} \Z \rightarrow G$ such that $h = \bar{h} \rho_F $,
where $\rho_F$ is the natural map from $\pmc{$\times$}\, \, \;_{i
\in \omega} \Z$ to $*_{i \in F} \Z$. Corresponding to the abelian
case this phenomenon was called the non-commutative Specker
Phenomenon. In this paper we show that Higman's result fails if
one passes from countable to uncountable. In particular, we show
that for non-trivial groups $G_{\alpha}$ $(\alpha \in \lambda)$
and uncountable cardinal $\lambda$ there are $2^{2^{\lambda}}$
homomorphisms from the complete free product of the
$G_{\alpha}$'s to the ring of integers.
\end{abstract}

\setcounter{section}{-1}

\section{Introduction}
Higman proved in 1952 \cite{H} that every free group $F$ is
non-commutatively slender, where slenderness means that any
homomorphism $h$ from the countable complete free product
$\pmc{$\times$}\, \, \;_{i \in \omega} \Z$ of the integers to $F$
depends only on finitely many coordinates. A similar result was
proved by Specker in 1950 \cite{S} for abelian groups. Specker
showed that any homomorphism from the countable product
$\Pi_{\omega} \Z$ to the integers depends only on finitely many
entries. These two phenomenons were called the commutative and
the non-commutative Specker Phenomenon. Eda extended Higman's
result in 1992 in \cite{E1} by showing that for any
non-commuatively slender group $S$, non-trivial groups
$G_{\alpha}$ ($\alpha \in I)$ and any homomorphism $h$ from the
free $\sigma$-product of the $G_{\alpha}'s$ to $\Z$ there exist a
finite subset $F$ of $I$ and a homomorphism $\bar{h} : *_{i \in
F} G_i \rightarrow S$ such that $h = \bar{h} \rho_F $, where
$\rho_F$ is the natural map from $\pmc{$\times$}\, \, \;_{i \in
I}^{\sigma} G_i$ to $*_{i \in F} G_i$. Motivated by this result
Eda asked the question \cite{E1}[Question 3.8] whether the
non-commutative Specker Phenomenon still holds if one passes from
countable to uncountable. In this paper we will answer this
question to the negative by constructing for a given uncountable
cardinal $\lambda$ and non-trivial groups $G_{\alpha}$ $(\alpha
\in \lambda)$, a homomorphisms $h$ from the complete free product
of the $G_{\alpha}$'s to $\Z$ for which the non-commutative
Specker Phenomenon fails. In particular, we will show that there
are $2^{2^{\lambda}}$ of these homomorphisms, hence the size of
the set of all homomorphisms from $\pmc{$\times$}\, \, \;_{\alpha
\in \lambda} G_{\alpha}$ to the integers is as large as possible.

\section{Basics and notations}

Let $I$ be an arbitrary indexset. For groups $G_i$ $(i \in I)$,
the free product is denoted by $*_{i \in I} G_i$ (see \cite{M}
for details on free products). If $J$ is a finite subset of $I$
then we write $J \Subset I$. For $X \subset Y \subset I$ let
$\rho_{XY} : *_{i \in Y}G_i \rightarrow *_{i \in X}G_i$ be the
canonical homomorphism. Then, the set $\{ *_{i \in X} G_i : X
\Subset I \}$ together with the homomorphisms $\rho_{XY}$ $(X
\subset Y \Subset I)$ form an inverse system and its inverse
limit $\underset{\leftarrow}{lim} (*_{i  \in X} G_i, \rho_{XY} :
X \subset Y \Subset I)$ is called the {\it unrestricted free
product} (see \cite{H}). Eda \cite{E1} introduced an infinite
version of free products and defined the {\it complete free
product} $\pmc{$\times$}\, \, \;_{i \in I}G_i$ of the groups
$G_i$ which is isomorphic to the subgroup $\bigcap_{F \Subset I}
\{*_{i \in F} G_i * \underset{\leftarrow}{lim}(*_{i \in X} G_i ,
\rho_{XY} : X \subset Y \Subset I) \}$ of the unrestricted free
product. To get familiar with the complete free product we recall
the definition of words of infinite length and some basic facts
about $\pmc{$\times$}\, \, \;_{i \in I}G_i$ from \cite{E1}.

\begin{definition}
Let $G_i$ $(i \in I)$ be non-trivial groups such that $G_i \cap
G_j =\{ e \}$ for $i \not= j \in I$. Elements of
$\bigcup\limits_{i \in I}G_i$ are called {\it letters}. A {\it
word} $W$ is a function $W: \Wb \rightarrow \bigcup\limits_{i \in
I}G_i$ such that $\Wb$ is a linearly ordered set and
$W^{-1}(G_i)$ is finite for any $i \in I$. In case the
cardinality of $\Wb$ is countable, we say that $W$ is a {\it
$\sigma$-word}. The class of all words is denoted by
${\cal{W}}(G_i : i \in I)$ (abbreviated by ${\cal{W}}$) and the
class of all $\sigma$-words is denoted by ${\cal{W}}^{\sigma}(G_i
: i \in I)$ (abbreviated by ${\cal{W}}^{\sigma}$).
\end{definition}

Two words $U$ and $V$ are said to be {\it isomorphic} $(U \cong
V)$ if there exists an isomorphism $\varphi : \Ub \rightarrow
\Vb$ as linearly ordered sets such that
$U(\alpha)=V(\varphi(\alpha))$ for all $\alpha \in \Ub$. It is
easily seen that ${\cal{W}}$ is a set and that for words of
finite length the above definition coincides with the usual
definition of words. For a subset $X \subset I$ the {\it
restricted word (or subword)} $W_X$ of $W$ is given by the
function $W_X: \WXb \rightarrow \bigcup\limits_{i \in X}G_i$,
where $\WXb = \{ \alpha \in \Wb : W(\alpha) \in \bigcup\limits_{i
\in X}G_i \}$ and $W_X(\alpha)=W(\alpha)$ for all $\alpha \in
\WXb$. Hence $W_X \in {\cal{W}}$. Now an equivalence relation is
defined on ${\cal{W}}$ by saying that two words $U$ and $V$ are
{\it equivalent} $(U \sim V)$ if $U_F = V_F$ for all $F \Subset
I$, where we regard $U_F$ and $V_F$ as elements of the free
product $*_{i \in F}G_i$. The equivalence class of a word $W$ is
denoted by $[W]$ and the composition of two words as well as the
inverse of a word are defined natural. Thus ${\cal{W}}/\sim$ $=
\{ [W] : W \in {\cal{W}} \}$ becomes a group.

\begin{definition}
The {\it complete free product} $\pmc{$\times$}\, \, \;_{i \in I}G_i$ is the
group ${\cal{W}}(G_i : i \in I)/\sim$. The {\it free $\sigma$-product}
$\pmc{$\times$}\, \, \;_{i \in I}^{\sigma}G_i$ is the group
${\cal{W}}^{\sigma}(G_i : i \in I)/\sim$, which is a subgroup of
$\pmc{$\times$}\, \, \;_{i \in I}G_i$. In case every $G_i$ is isomorphic to
$G$, we abbreviate $\pmc{$\times$}\, \, \;_{i \in I}G_i$ by
$\pmc{$\times$}\, \, \;_IG$ and similarly for free $\sigma$-products.
\end{definition}

Obviously, $\pmc{$\times$}\, \, \;_{i \in I}G_i$ and
$\pmc{$\times$}\, \, \;_{i \in I}^{\sigma}G_i$ are isomorphic to
$*_{i \in I}G_i$ if $I$ is finite. By \cite[Proposition 1.8]{E1}
the complete free product $\pmc{$\times$}\, \, \;_{i \in I}G_i$
is isomorphic to the subgroup $\bigcap_{F \Subset I} \{*_{i \in
F} G_i * \underset{\leftarrow}{lim}(*_{i \in X} G_i , \rho_{XY} :
X \subset Y \Subset I) \}$ of the unrestricted free product.
Moreover, Eda proved in \cite{E1} that each equivalence class
$[W]$ is determined uniquely by a reduced word. A word $W \in
{\cal{W}}(G_i : i \in I)$ is called {\it reduced}, if $W \cong
UXV$ implies $[X] \not= e$ for any non-empty word $X$, where $e$
is the identity, and for any neighboring elements $\alpha$ and
$\beta$ of $\Wb$ it never occurs that $W(\alpha)$ and $W(\beta)$
belong to the same $G_i$.

\begin{lemma}[Eda, \cite{E1} ]
\label{reduced} For any word $W \in {\cal{W}}(G_i : i \in I)$,
there exists a reduced word $V \in {\cal{W}}(G_i : i \in I)$ such
that $[W]=[V]$ and $V$ is unique up to isomorphism.
\end{lemma}

Furthermore, Eda showed in \cite{E1} the following lemma where a
word $W \in {\cal{W}}(G_i : i \in I)$ is called {\it
quasi-reduced} if the reduced word of $W$ can be obtained by
multiplying neighboring elements without cancelation.

\begin{lemma}[Eda, \cite{E1}]
\label{quasireduced} For any two reduced words $W, V \in
{\cal{W}}(G_i : i \in I)$ there exist reduced words $V_1,W_1,M
\in {\cal{W}}(G_i : i \in I)$ such that $W \cong W_1M$, $V \cong
M^{-1}V_1$ and $W_1V_1$ is quasi-reduced.
\end{lemma}

We would like to remark that a free $\sigma$-product
$\pmc{$\times$}\, \, \;_{i \in I}^{\sigma}\Z_i$ is isomorphic to
the fundamental group and the free complete product
$\pmc{$\times$}\, \, \;_{i \in I}\Z_i$ is isomorphic to the big
fundamental group of the Hawaiian earring with $I$-many circles
(see \cite{CC}). Hence free complete products are also of
topological interest.

\section{The uncountable Specker Phenomenon}
In 1950, E. Specker \cite{S} proved that for any homomorphism $h$
from the countable direct product $\mathbb{Z}^{\omega}$ to the
ring of integers $\mathbb{Z}$, there exist a finite subset $F$ of
$\omega$ and a homomorphism $\bar{h} : \Z^F \rightarrow
\mathbb{Z}$ satisfying $h = \bar{h}\rho_F$, where $\rho_F:
\Z^{\omega} \rightarrow \Z^F$ is the canonical projection. This
phenomenon is called the Specker Phenomenon and it can easily be
seen that Specker's result still holds if one replaces
$h:\mathbb{Z}^{\omega} \rightarrow \mathbb{Z}$ by $g: \Z^{\omega}
\rightarrow G$, where $G$ is any free abelian group. For
generalizations to products of larger cardinalities and the
resulting definition of slenderness for abelian groups we refer
to \cite{EM} or \cite{F1}. In \cite{E2} Eda introduced a
non-commutative version of slenderness.

\begin{definition}
A group $G$ is {\it non-commutatively slender} if for any
homomorphism $h: \pmc{$\times$}\, \, \;_{\N} \Z \rightarrow G$
there exists a natural number $n$ such that $h(\pmc{$\times$}\,
\, \;_{\N \backslash \{1,\cdots,n \} } \Z) = \{ e \}$.
\end{definition}

Eda proved that every non-commutatively slender group is
torsion-free and that in the abelian case non-commuative
slenderness is equivalent to the commutative slenderness (see
\cite[Theorem 3.3. and Corollary 3.4.]{E1}). Moreover, he proved
that non-commutative slender groups behave nice in the following
sense for non-trivial groups $G_i$ $(i \in I)$.

\begin{proposition}[Eda, \cite{E1}]
Let $S$ be a non-commutative slender group and $h :
\pmc{$\times$}\, \, \;_{i \in I}^{\sigma} G_i \rightarrow S$ be a
homomorphism. Then, there exist a finite subset $F$ of $I$ and a
homomorphism $\bar{h} : *_{i \in F} G_i \rightarrow S$ such that
$h = \bar{h} \rho_F $, where $\rho_F$ is the natural map from
$\pmc{$\times$}\, \, \;_{i \in I}^{\sigma} G_i$ to $*_{i \in F}
G_i$.
\end{proposition}

Moreover, if $S_j$ $(j \in J)$ are non-commutatively slender
groups then also the restricted direct product and the free
product of the $S_j$'s are non-commutatively slender (see
\cite[Theorem 3.6.]{E1}). The first fundamental result on the
class of non-commutatively slender groups was already obtained by
Higman in \cite{H} where he proved the following theorem.

\begin{theorem}[Higman, \cite{H}]
Every free group is non-commutatively slender.
\end{theorem}

In contrast to Higman's result we will show that if one replaces
countable by uncountable then the non-commutative Specker
Phenomenon fails. In particular we show that there are
$2^{2^{\lambda}}$ homomorphisms from the complete free product
$\pmc{$\times$}\, \, \;_{\alpha \in \lambda}G_{\alpha}$ to the
ring of integers if $\lambda$ is uncountable and the
$G_{\alpha}$'s are non-trivial groups. For the convenience of the
reader we first construct one homomorphism for which the Specker
Phenomenon fails and then modify the construction to obtain our
main result.

\begin{theorem}
\label{main1} Let $\lambda$ be any uncountable cardinal and
$G_{\alpha}$ $(\alpha \in \lambda)$ non-trivial groups. For each
$\kappa \leq \lambda$ regular, uncountable there exists a
homomorphism $\varphi_{\kappa} : \pmc{$\times$}\, \, \;_{\alpha
\in \lambda}G_{\alpha} \rightarrow \Z$ for which the Specker
Phenomenon fails.
\end{theorem}

\proof Let $G_{\alpha}$ $(\alpha < \lambda)$ be a collection of
non-trivial groups and choose $e_{\alpha}\not= g_{\alpha} \in
G_{\alpha}$, where $e_{\alpha}$ is the identity of $G_{\alpha}$
$(\alpha < \lambda)$. We define the following words $M_{\kappa}
\in \pmc{$\times$}\, \, \;_{\alpha \in \lambda}G_{\alpha}$ for any
regular, uncountable cardinal $\kappa < \lambda$.
\[ M_{\kappa} : (\kappa, <) \longrightarrow \cup_{\alpha < \lambda}
G_{\alpha} \text{ via } \beta \mapsto g_{\beta} \]
where $<$ is the natural ordering of $\lambda$. Note that
$M_{\kappa}$ is a word of uncountable cofinality since $\kappa$
is regular and uncountable. For $\beta < \kappa$ we let
$M_{\kappa, \beta}$ be the subword
$M_{\kappa}\restriction_{[\beta, \kappa)}$ of $M_{\kappa}$. Now
let $X$ be any reduced word in $ \pmc{$\times$}\, \, \;_{\alpha
\in \lambda}G_{\alpha}$ and recall that a subset $J \subseteq
(\lambda, <)$ is called convex if $x < y < z$ and $x,z \in J$
implies $y \in J$. We put
\[ Occ^+_{\kappa}(X):= \{ J \subseteq (\lambda,<) : J \text{ is convex and }
X \restriction_J \cong M_{\kappa, \beta} \text{ for some } \beta < \kappa
\}. \]
Thus $Occ^+_{\kappa}(X)$ counts the occurencies of end segments of
$M_{\kappa}$ in $X$. Similarly we let
\[ Occ^-_{\kappa}(X):= \{ J \subseteq (\lambda,<) : J \text{ is convex and }
X\restriction_J \cong M^{-1}_{\kappa, \beta} \text{ for some } \beta <
\kappa \}. \]
In order to avoid counting subsets of $(\lambda,<)$ several times
we define an equivalence relation on $Occ^+_{\kappa}(X)$ and
$Occ^-_{\kappa}(X)$ in the following way. Two convex subsets
$J_1,J_2$ of $(\lambda,<)$ are said to be equivalent if they have
a common end segment, i.e. $J_1 \sim_{\kappa} J_2$ if there exist
$j_i \in J_i$ such that $X \restriction_{S_1} \cong X
\restriction_{S_2}$, where $S_i =\{ j \in J_i : j \geq j_i \}$.
First we prove that two subsets $J_1,J_2 \in Occ^+_{\kappa}(X)$
are either disjoint or equivalent. Therefore assume that $J_1,J_2
\in Occ^+_{\kappa}(X)$ are not disjoint, hence there exists $j^*
\in J_1 \cap J_2$. We let $h_i : M_{\kappa,\beta_1}
\longrightarrow X \restriction_{J_i}$ be isomorphisms for some
$\beta_i < \lambda$ $(i=1,2)$. Thus we can find $\gamma_i \geq
\beta_i$ such that $h_i(\gamma_i)=j^*$ and therefore
$X(j^*)=g_{\gamma_i}$ for $i=1,2$. Hence $\gamma_1=\gamma_2$ and
by transfinite induction we conclude $X \restriction_{T_1} \cong X
\restriction_{T_2}$, where $T_i = \{ j \in J_i : j \geq j^* \}$.
Note that $h_i$ is an isomorphism of linearly ordered sets, hence
$h_i$ commutes with limits and the successor-function. Similarly
two subsets $J_1,J_2 \in Occ^-_{\kappa}(X)$ are either disjoint or
equivalent. Next we will show that the sets
$Occ^+_{\kappa}(X)/\sim_{\kappa}$ and
$Occ^-_{\kappa}(X)/\sim_{\kappa}$ are finite. Therefore assume
that there exist infinitely many pairwise non-equivalent $J_n \in
Occ^+_{\kappa}(X)$ $(n \in \omega)$. Hence $J_n$ and $J_m$ are
disjoint for $n \not= m$. We let $X\restriction_{J_n} \cong
M_{\kappa,\beta_n}$ for some $\beta_n < \kappa$ and $n \in
\omega$. Then $\beta = \cup_{n \in \omega}\beta_n$ is strictly
less than $\kappa$ since $\kappa$ is regular and uncountable,
hence $cf(\kappa) > \aleph_0$. Since $\beta \in [\beta_n,
\kappa)$ for all $n \in \omega$ we can find $j_n \in J_n$ such
that
\[ X(j_n)=M_{\kappa, \beta_n}(\beta)=M_{\kappa, \beta}(\beta). \]
for $n \in \omega$. But all $J_n$ are pairwise disjoint and
therefore $X^{-1}(G_{\beta})$ is infinite which is a
contradiction. Thus $Occ^+_{\kappa}(X)/\sim_{\kappa}$ and
similarly $Occ^-_{\kappa}(X)/\sim_{\kappa}$ are finite sets. We
now define $\varphi_{\kappa}: \pmc{$\times$}\, \, \;_{\alpha \in
\lambda}G_{\alpha} \longrightarrow \Z$ as follows:
\[ X \mapsto |Occ^+_{\kappa}(V)/\sim_{\kappa}| -
|Occ^-_{\kappa}(V)/\sim_{\kappa} | \]
where $V$ is the reduced word corresponding to $X$.\\ Note that
$\varphi_{\kappa}$ is well-defined by Lemma \ref{reduced}.
Moreover, by definition $\varphi_{\kappa}(X^{-1}) = -
\varphi_{\kappa}(X)$ and obviously the Specker Phenomenon fails
for $\varphi_{\kappa}$. All we have to show is that
$\varphi_{\kappa}$ is a homomorphism. Therefore let $X$ and $Y$
be reduced words. By Lemma \ref{quasireduced} there exist reduced
words $X_1,Y_1$ and $M$ such that $X \cong X_1M$ and $Y \cong
M^{-1}Y_1$ and $X_1Y_1$ is quasi-reduced. Now it is easy to check
that $ \varphi_{\kappa}(XY) = \varphi_{\kappa}(X_1Y_1)$ by
definition and hence
\[ \varphi_{\kappa}(XY) = \varphi_{\kappa}(X_1Y_1)= \varphi_{\kappa}(X_1) +
\varphi_{\kappa}(Y_1) = \varphi_{\kappa}(X) + \varphi_{\kappa}(Y), \]
since $X_1Y_1$ is quasi-reduced.
\qed

We would like to remark that the uncountable cofinality of
$\lambda$ in Theorem \ref{main1} is essential and can not be
avoided by Higman's theorem. Modifying the proof of Theorem
\ref{main1} we obtain

\begin{theorem}
\label{main2} Let $\lambda$ be an uncountable cardinal and
$G_{\alpha}$ $( \alpha \in \lambda)$ be non-trivial groups. Then
there are $2^{2^{\lambda}}$ homomorphisms from the complete free
product of the $G_{\alpha}$'s to the ring of integers.
\end{theorem}

\proof Let $M_{\alpha}$ be a reduced word in $\pmc{$\times$}\, \,
\;_{\alpha \in \lambda}G_{\alpha}$ of uncountable cofinality
$\lambda$, i.e. $\bar{M_{\alpha}} = (\lambda,<)$. Recall from the
proof of Theorem \ref{main1} that by $M_{\alpha,\beta}$ we mean
the subword $M_{\alpha} \restriction_{[\beta,\lambda)}$ for $\beta
\in \lambda$. Assume that we have a family of such words
$M_{\alpha}$ $(\alpha \in 2^{\lambda})$ satisfying the following
condition for a convex subset $J \subseteq \lambda$ and a reduced
word $X$:
\begin{equation}
\tag{$*$} X \restriction_J \cong M_{\alpha,\beta} \text{ for }
\beta \in \lambda \Longrightarrow X \restriction_J \not\cong
M_{\gamma, \delta} \text{ for all } \alpha \not= \gamma \in
2^{\lambda}, \delta \in \lambda
\end{equation}
Then it is well-known that we can choose $2^{2^{\lambda}}$ almost
disjoint families $F_\alpha = \{ M_{\beta} : \beta \in I_{\alpha}
\}$ such that $I_{\alpha}$ has size $\lambda$. We now define for
$\alpha \in 2^{2^{\lambda}}$
\[ Occ^+_{\alpha}(X):= \{ J \subseteq (\lambda,<) : J \text{ convex, }
 X \restriction_J \cong M_{\beta, \gamma} \text{ for some } \beta \in
I_{\alpha}, \gamma < \lambda \} \]
and similarly
\[ Occ^-_{\alpha}(X):= \{ J \subseteq (\lambda,<) : J \text{ convex, }
X\restriction_J \cong M^{-1}_{\beta, \gamma} \text{ for some }
\beta \in I_{\alpha}, \gamma < \lambda \}. \] As in the proof of
Theorem \ref{main1} we can see that the sets
$Occ^+_{\alpha}(X)/\sim_{\lambda}$ and
$Occ^-_{\alpha}(X)/\sim_{\lambda}$ are finite for any reduced word
$X$ and $\alpha \in 2^{2^{\lambda}}$. Moreover, the maps
$\varphi_{\alpha}: \pmc{$\times$}\, \, \;_{\beta \in
\lambda}G_{\beta} \longrightarrow \Z$ defined by
\[ X \mapsto |Occ^+_{\alpha}(V)/\sim_{\lambda}| -
|Occ^-_{\alpha}(V)/\sim_{\lambda} | \]
where $V$ is the reduced word corresponding to $X$, are
well-defined homomorphisms and $\varphi_{\alpha} \not=
\varphi_{\beta}$ for $\alpha, \beta \in 2^{2^{\lambda}}$ since
the families $F_{\delta}$ are almost disjoint and satisfy
condition $(*)$. Hence the size of all homomorphisms from the
complete free product of the $G_{\alpha}$'s to the integers is
$2^{2^{\lambda}}$ as claimed.
It remains to show the existence of the words $M_{\alpha}$ satisfying
$(*)$.\\
We start with any partition of $\lambda$ into two sets, i.e. with
a function $g: \lambda \rightarrow \{ 0,1 \}$. Note that there
are $2^{\lambda}$ of those functions. Moreover, we choose
elements $e_{\alpha} \not= h_{\alpha} \in G_{\alpha}$ $(\alpha \in
\lambda)$ and define the word $M_g^{\prime} \in \pmc{$\times$}\,
\, \;_{\alpha \in \lambda}G_{\alpha}$ by
\[ M_g^{\prime} (\beta) = h_{\beta + 2g(\beta)} \]
Then $M_g^{\prime}$ is a reduced word and we let $M_g$ be the
composition of $M_g^{\prime}$ with itself $\omega_1$ times. Then
$M_g$ is still reduced and for different $g,g^{\prime} : \lambda
\rightarrow \{ 0,1 \}$ condition $(*)$ is satisfied for $M_g$ and
$M_{g^{\prime}}$. Thus the family
\[ F = \{ M_g : g: \lambda \rightarrow \{ 0,1 \} \} \]
is a family of reduced words of size $2^{\lambda}$ satisfying
$(*)$ as desired. \qed

Moreover, the authors would like to mention that modifying the
proof of Theorem \ref{main1} Conner and Eda proved a more general
result in \cite{CE}

\goodbreak

\bigskip

\noindent
Saharon Shelah\\
Department of Mathematics\\
Hebrew University\\
Jerusalem, Israel\\
and Rutgers University\\
Newbrunswick, NJ U.S.A.\\
e-Mail: Shelah@@math.huji.ae.il\\

\noindent
Lutz Str\"ungmann\\
Fachbereich 6, Mathematik\\
University of Essen\\
45117 Essen\\
Germany\\
e-Mail: lutz.struengmann@@uni-essen.de\\

\begin{thebibliography}{99}
\bibitem[CC]{CC}{\bf J. W. Cannon and G. Conner}, The big fundamental group,
big Hawaiian earrings, and the big free groups, to appear in Top. Appl.
Math.

\bibitem[CE]{CE}{\bf G. Conner and K. Eda}, Free subgroups of complete free
products, preprint (2000).

\bibitem[E1]{E1}{\bf K. Eda}, Free $\sigma$-products and noncommutatively
slender groups, {\it J. of Algebra} {\bf 148} (1992), 243--263.

\bibitem[E2]{E2}{\bf K. Eda}, The non-commutative Specker Phenomenon, {\it
J. of Algebra} {\bf 204} (1998), 95--107.

\bibitem[EM]{EM}{\bf P.C. Eklof and A. Mekler}, Almost Free Modules, Set-
Theoretic Methods, Amsterdam, New York, North-Holland, Math. Library.

\bibitem[F1]{F1}{\bf L. Fuchs}, Infinite Abelian Groups - Volume I and II,
Academic Press, New York-London (1970 and 1973).

\bibitem[F2]{F2}{\bf L. Fuchs}, Abelian Groups, Hungarian Academy of
Science, Budapest (1958).

\bibitem[H]{H}{\bf G. Higman}, Unrestricted free products and varieties of
topological groups, {\it J. London Math. Soc.} {\bf 27} (1952), 73--81.

\bibitem[M]{M}{\bf W.S. Massey}, Algebraic Topology: An Introduction, {\it
Graduate Texts in Mathematics}, Springer Verlag New York-Heidelberg-Berlin
(1967).


\bibitem[S]{S}{\bf E. Specker}, Additive Gruppen von Folgen ganzer Zahlen,
{\it Portugal. Math.} {\bf 9} (1950), 131--140.

\end{thebibliography}
\end{document}